\documentclass[11pt]{amsart}

\newcommand{\w}{\omega}
\newcommand{\e}{\varepsilon}

\newtheorem{theorem}{Theorem}

\title[Non-compact operators into $c_0$]{Constructing non-compact operators into $c_0$}
\author{Iryna Banakh and Taras Banakh}
\address{Department of Functional Analysis, Ya.Pidstryhach Institute for Applied Problems of Mechanics and Mathematics, Naukova 3b, Lviv, Ukraine}
\email{ibanakh@yahoo.com}
\address{Instytut Matematyki, Uniwersytet Humanistyczno-Przyrodniczy Jana Kochanowskiego, Kielce, Poland, and Department of Mathematics, Ivan Franko National University of Lviv, Universytetska 1, 79000, Lviv, Ukraine}
\email{t.o.banakh@gmail.com}

\begin{document}
\begin{abstract} We prove that for each dense non-compact linear operator $S:X\to Y$ between Banach spaces there is a linear operator $T:Y\to c_0$ such that the operator $TS:X\to c_0$ is not compact. This generalizes the Josefson-Nissenzweig Theorem.
\end{abstract}
\maketitle

By the Josefson-Nissenzweig Theorem \cite{Jos}, \cite{Nis} (see also \cite{Beh}, \cite[XII]{Di}, and \cite[3.27]{HMVZ}), for each infinite-dimensional Banach space $Y$ 
the weak$^*$ convergence and norm convergence in the dual Banach space $Y^*$ are distinct. This allows us to find a sequence $(y_n^*)_{n\in\w}$ of norm-one functionals in $Y^*$ that converges to zero in the weak$^*$ topology. 
Such functionals determine a non-compact operator $T:Y\to c_0$ that assigns to each $y\in Y$ the vanishing sequence $(y_n^*(y))_{n\in\w}\in c_0$. Thus each infinite-dimensional Banach space $Y$ admits a non-compact  operator $T:Y\to c_0$ into the Banach space $c_0$. 

The following theorem (which is a crucial ingredient in the topological classification \cite{BC} of closed convex sets in Fr\'echet spaces) says a bit more:

\begin{theorem}\label{main} For any dense non-compact operator $S:X\to Y$ between Banach spaces
there is an operator $T:Y\to c_0$ such that the composition $TS:X\to c_0$ is non-compact.
\end{theorem}

By an {\em operator} we understand a linear continuous operator. An operator $T:X\to Y$ is {\em dense} if $T(X)$ is dense in $Y$.

The proof of Theorem~\ref{main} uses the famous Rosenthal $\ell_1$ Theorem \cite{Ros} (see also \cite[XI]{Di} and \cite{Beh}) saying that any bounded sequence in a Banach space $X$ contains a subsequence which is either weakly Cauchy or $\ell_1$-basic. 

A sequence $(x_n)_{n\in\w}$ in a Banach space $(X,\|\cdot\|)$ is called {\em $\ell_1$-basic} if there are constants $0<c\le C<\infty$ such that for each real sequence $(\alpha_n)_{n\in\w}\in \ell_1$ we get
$$c\sum_{n\in\w}|\alpha_n|\le \Big\|\sum_{n\in\w}\alpha_nx_n\Big\|\le C\sum_{n\in\w}|\alpha_n|.$$ 

\begin{proof}[Proof of Theorem~\ref{main}] Assume that $S:X\to Y$ is a dense non-compact operator. Let $(e_n)_{n\in\w}$ be the standard Schauder basis of the Banach space $c_0$ and $(e_n^*)_{n\in\w}$ is the dual basis in the dual space $c_0^*=\ell_1$. 
To construct the operator $T:Y\to c_0$ with non-compact $TS$, we shall consider three cases.

1. First we assume that the following condition holds:
\begin{itemize}
\item[(i)] there is an $\ell_1$-basic sequence $(y_n^*)_{n\in\w}$ in $Y^*$ such that the sequence $(S^*y^*_n)_{n\in\w}$ is $\ell_1$-basic and weak$^*$ null in $X^*$.
\end{itemize}

In this case we define the operator $T:Y\to c_0$ by $T:y\mapsto(y^*_n(y))_{n\in\w}$. Observe that the dual operator $T^*:c_0^*\to Y^*$ maps the $n$-th coordinate functional $e^*_n\in c_0^*$ onto $y^*_n$. Consequently, the sequence
$$\big(S^*y^*_n\big)_{n\in\w}=\big((TS)^*e^*_n\big)_{n\in\w},$$
being $\ell_1$-basic, is not totally bounded in $Y^*$, which implies that the dual operator $(TS)^*:c_0^*\to X^*$ is not compact. By the Schauder Theorem 
\cite[7.7.]{FHHMPZ}, the operator $TS:X\to c_0$ also is not compact.
\smallskip

2. Assume that the condition (i) does not hold but 
\begin{itemize}
\item[(ii)] there is an $\ell_1$-basic sequence $(y_n^*)_{n\in\w}$ in $Y^*$ whose image $(S^*y^*_n)_{n\in\w}$ is $\ell_1$-basic in $X^*$.
\end{itemize}
In this case, by \cite[Exercise 3(i)]{Di} the condition (ii) combined with the negation of (i) imply the existence of an $\ell_1$-basic sequence $(x_n)_{n\in\w}$ in $X$ whose image $(Sx_n)_{n\in\w}$ is an $\ell_1$-basic sequence in $Y$.
Arguing as in the proof of Josefson-Nissenzweig Theorem \cite[p.223]{Di}, we can construct a bounded linear operator $T:Y\to c_0$ such that $TS(x_n)=e_n\in c_0$ for all $n\in\w$. Since the operator $TS$ is not compact, we are done.
\smallskip

3. Assume that (ii) does not hold. Since the operator $S$ is not compact, its dual $S^*:Y^*\to X^*$ is not compact too, see \cite[7.7]{FHHMPZ}.
This means that the image $S^*(B^*)$ of the closed unit ball $B^*\subset Y^*$ is not totally bounded in $X^*$. Consequently, the dual ball $B^*$ contains a sequence $(y_n^*)_{n\in\w}$ whose image $(S^*y^*_n)_{n\in\w}$ is $\e$-separated for some $\e>0$. The latter means that $\|S^*(y_n^*-y^*_m)\|\ge\e$ for all $n\ne m$.

By the Rosenthal $\ell_1$ Theorem, $(S^*y_n^*)_{n\in\w}$ contains a subsequence which is either weak Cauchy or $\ell_1$-basic. We lose no generality assuming that the entire sequence $(S^*y_n^*)_{n\in\w}$ is either weak Cauchy or $\ell_1$-basic. 
\smallskip

3a. First we assume that the sequence $(S^*y^*_n)_{n\in\w}$ is weak Cauchy. 
Then it is weak$^*$ Cauchy and being a subset of the weakly$^*$ compact set $S^*(B^*)$ weakly$^*$ converges to some point $x^*_\infty\in S^*(B^*)$. Fix any point $y^*_\infty\in B^*$ with $S^*(y^*_\infty)=x^*_\infty$. 
The density of the operator $S:X\to Y$ implies the injectivity of the dual operator $S^*:Y^*\to X^*$. The weak$^*$ compactness of the closed unit ball $B^*\subset Y^*$ guarantees that $S^*|B^*:B^*\to X^*$ is a homeomorphic embedding for the weak$^*$ topologies on $B^*$ and $X^*$. Now we see that the weak$^*$ convergence of the sequence $(S^*y^*_n)_{n\in\w}$ to $S^*y^*_\infty$ implies the weak$^*$ convergence of the sequence $(y^*_n-y^*_\infty)_{n\in\w}$ to zero.

Then the bounded operator  $T:Y\to c_0$, $T:y\mapsto\big((y^*_n-y^*_\infty)(y)\big)_{n\in\w}$, is well-defined. Since the set $\{(TS)^{*}(e^*_n)\}_{n\in\w}=\{S^*(y^*_n-y^*_\infty)\}_{n\in\w}$ is $\e$-separated, the operator $(TS)^*:c_0^*\to X^*$ is not compact and hence $TS:X\to c_0$ is not compact too.
\smallskip

3b. Finally, assume that $(S^*y_n^*)_{n\in\w}$ is an $\ell_1$-basic sequence in $X^*$.  By Proposition 5.10 \cite{FHHMPZ} (the lifting property of  $\ell_1$), the sequence $(y_n^*)_{n\in\w}$ is $\ell_1$-basic in $Y^*$, which contradicts our assumption that the condition (ii) fails.
\end{proof}


\begin{thebibliography}{}

\bibitem{BC} T.~Banakh, R.~Cauty, {\em Topological classification of closed convex sets in Fr\'echet spaces}, preprint.

\bibitem{Beh} E.~Behrends, {\em New proofs of Rosenthal's $l\sp 1$-theorem and the Josefson-Nissenzweig theorem}, Bull. Polish Acad. Sci. Math. {\bf 43}:4 (1995), 283--295 (1996).

\bibitem{HMVZ} P.~Hajek, V.S.~Montesinos,J.~Vanderwerff, V.~Zizler, Biorthogonal systems in Banach spaces, Springer, NY, 2008. 

\bibitem{FHHMPZ} M.~Fabian, P.~Habala, P.~Hajek, V.S.~Montesinos, J.~Pelant, V.~Zizler, Functional analysis and infinite-dimensional geometry, Springer, NY, 2001.

\bibitem{Di} J.~Diestel, Sequences and series in Banach spaces, Springer, NY, 1984. 

\bibitem{Jos} B.~Josefson, {\em Weak sequential convergence in the dual of a Banach space does not imply norm convergence}, Ark. Mat. {\bf 13} (1975), 79--89.

\bibitem{Nis} A.~Nissenzweig, {\em $W\sp{\ast} $ sequential convergence}, Israel J. Math. {\bf 22}:3-4 (1975), 266--272.

\bibitem{Ros} H.P.~Rosenthal, {\em A characterization of Banach spaces containing $\ell_1$,} Proc. Natl. Acad. Sci. USA {\bf 71} (1974), 2411-2413.

\end{thebibliography}
\end{document}